\documentclass[twoside,10pt, leqno]{article}
\usepackage[spanish, mexico]{babel}
\usepackage[utf8]{inputenc}
\usepackage{multicol}
\usepackage{graphicx}
\usepackage{amssymb}
\usepackage{amscd, amsfonts, amsmath, amsthm, amssymb}
\usepackage{enumerate}
\usepackage{color}
\usepackage{mathtools}
\usepackage{tabularx}
\usepackage{longtable}
\usepackage{enumerate}
\usepackage{booktabs}
\usepackage{comment}
\usepackage[table]{xcolor}
\usepackage{cases}
\usepackage{pdflscape}
\usepackage{xspace}
\usepackage{bbm}
\usepackage[caption=false]{subfig}
\usepackage{nicefrac}
\usepackage{graphicx}
\usepackage{cite}  
\usepackage{nccfoots}
\usepackage{cancel}
\usepackage[small]{caption2}
\usepackage{multirow}
\usepackage{enumerate}

\usepackage{tabulary}
\usepackage{array}
\usepackage{comment}
\usepackage{diagbox}
\usepackage{hyperref}
\usepackage{array}

\input xy
\xyoption{all}

\theoremstyle{plain}
\newtheorem{teo}{Theorem}[section]
\newtheorem{cor}[teo]{Corollary}
\newtheorem{lema}[teo]{Lemma}
\newtheorem{prop}[teo]{Proposition}
\theoremstyle{definition}
\newtheorem{defi}[teo]{Definition}
\newtheorem{obs}[teo]{Remark}
\newtheorem{ejem}[teo]{Example}

\newcommand{\nphantom}[1]{\sbox0{#1}\hspace{-\the\wd0}}

\numberwithin{equation}{section}

\begin{document}
\thispagestyle{plain}

\bigskip

\begin{center}

{\bf{\Large A second note on homological systems}}\\

\smallskip

{\large {$^a$Jes\'us Efr\'en P\'erez Terrazas, $^b$Luis Fernando Tzec Poot}} 

$^{a,b}$ Facultad de Matemáticas, Universidad Autónoma de Yucat\'an\\

\medskip

$^a$jperezt@correo.uady.mx

\end{center}

\renewcommand{\abstractname}{Abstract}
\begin{abstract}

Let $\left( \Delta; \Omega , \leq \right)$ be a b-homological system and $\widetilde{\cal F} \left( \Delta \right)$ the category of the extended $\Delta -$filtered modules. Here there is a proof that $\widetilde{\cal F} \left( \Delta \right)$ is closed under direct summands.  
\end{abstract}

$\underline{ \ \ \ \ \ \ \ \ \ \ \  \ \ \ \ \ \ \ \ \ \ \ \ \ \ \ \ \ }$

{\footnotesize $Keywords$ $and$ $phrases:$ homological system, preorder, extended filtered modules.}

{\footnotesize 2020 $Mathematics$ $Subject$ $Classification:$ $16G99,$ $16E99$}

$\overline{ \ \ \ \ \ \ \ \ \ \ \  \ \ \ \ \ \ \ \ \ \ \ \ \ \ \ \ \ }$

\section {Introduction}

In the study of the representation type of the category of $\Delta-$filtered modules of an homological system, we can use the notion of generic module (see \cite{BPS}), and then it becomes natural to introduce the notion of the category of extended $\Delta -$filtered modules, and to understand some of its properties.

Here we give a proof of the closure under direct summands of the second category, adapting the arguments of \cite{P} to prove the closure under direct summands of the first.

\bigskip

\section{Pre-orders and filtrations}

Let $C$ be a not empty set. The relation $\leq$ on $C$ is a \emph{pre-order} if it is reflexive and transitive. In this case $\left( C, \leq \right)$ is a \emph{pre-ordered set.}

The following claims are easy to prove and also are not hard to find in internet:

\begin{lema}\label{growing} $\,$

\begin{enumerate}

\item There is an equivalence relation on $C$ given by $c_1 \sim c_2$ iff $c_1 \leq c_2$ and $c_2 \leq c_1.$ Let $\pi : C \longrightarrow C / \sim$ be the canonical induced suprajective function, and let us denote $\pi (c) = \overline{c}.$ Then $C / \sim$ has a canonical partial order given by $\overline{c_1} \leq \overline{c_2}$ iff $c_1 \leq c_2 .$ (See section 2.1 of \cite{MSX})

\item By the order-extension principle, \cite{S}, there is a total order (a linearization) $\preceq$ on $C / \sim$ that extends $\leq .$ 

\end{enumerate}

\end{lema}

\smallskip

With respect to Lemma \ref{growing} let $\left( C , \leq \right)$ denote the pre-ordered set, $\left( C / \sim , \leq \right)$ the partially ordered set and $\left( C / \sim  , \preceq \right)$ the totally ordered set.

\begin{lema} \label{related} Consider the setting of Lemma \ref{growing} and previous notation. 

\begin{enumerate}

\item If $c_1 \leq c_2$ then $\overline{c_1} \leq \overline{c_2}$ and $\overline{c_1} \preceq \overline{c_2}.$ 

\item By the proof ot the order-extension principle, if $\overline{c_1}$ is not related to $\overline{c_2}$ in the partially ordered set, then there is at least one total order that extends $\leq$ such that $\overline{c_1} \preceq \overline{c_2},$ and at least one total order that extends $\leq$ such that $\overline{c_2} \preceq \overline{c_1}.$

\end{enumerate}

\end{lema}

\smallskip

\begin{ejem} \label{naturals} Consider the set of natural numbers $\mathbb{N}$ and the partial order given by $n \leq m$ iff \lq\lq $n$ divides $m$\rq\rq.

In this partial order 1 is the minimal element.

There are infinite ways to extend the previous partial order to a total order, and here we want to consider some lexicographic orders: 

1.- Let $q (n)$ be the number of prime factors, including repetitions, in a factorization in primes of $n.$

2.- Choose an arbitrary total order $\preceq _1$ for those naturals with $q$ equal to one, some arbitrary total order $\preceq _2$ for those naturals with $q$ equal to 2, etc.  

3.- Given $n,m \in \mathbb{N}$ if $q(n) < q(m)$ then $n \preceq m,$ and if $q (n) = q (m) = q$ then $n \preceq m$ iff $n \preceq _q m.$ 

\end{ejem}

\smallskip

\begin{ejem} \label{inverter} Consider $n, m \in \mathbb{N}$ such that $n$ does not divide $m$ and $m$ does not divide $n.$ If $m$ has more prime factors that $n$ then in any lexicographic order of Example \ref{naturals} we get $n \prec m.$ If $q (m) = q(n) = a$ then we can choose $\preceq _a$ such that $m \prec n$ or $n \prec m .$

Now an idea to have $q (m) < q(n)$ and get $n \prec m .$

First observe that given an injective function $\iota : A \longrightarrow \mathbb{N}$ the total order $\preceq$ induces a total order on $A.$ 

Let $d = \left( n, m \right)$ be the maximal common divisor of $n$ and $m,$ and $n = n' d$ and $m = m' d,$ where $\left( n', m' \right) = 1.$ 

Let $i : \mathbb{N} \longrightarrow \mathbb{N}$ be the function given by $i (a) = a$ if $m'$ does not divide $a$ and $i (a) = n' a$ if $m' | a,$ and choose for the codomain a total order as in Example \ref{naturals}. It is easy to verify that $i$ is injective and also that $i$ induces a total order $\preceq$ on $\mathbb{N}$ such that $n' \prec m'$ and $n \prec m,$ and such that it also extends the partial order \lq\lq divides to\rq\rq.

\end{ejem}

\smallskip

The next definition is almost the same given in \cite{MSX}. 

$R$ will denote a ring.

\begin{defi} 
A b-homological system $\left(\Delta ; \Omega , \leq
\right),$ consists of the following:

HS1: A pre-ordered set $\left( \Omega , \leq \right) .$

HS2: The set $\Delta = \left\{ \Delta _{\omega} \right\}_{\omega \in
\Omega},$ where $\Delta _{\omega}$ is a finitely generated $R-$module with local ring of endomorphisms for any $\omega \in \Omega .$ 
Also $\omega \neq \omega '$ implies $\Delta _{\omega} \ncong
\Delta _{\omega '}.$

HS3: If ${\rm Hom}_{\Lambda} \left( \Delta _{\omega} , \Delta
_{\omega '} \right) \neq 0$ then $\omega \leq \omega ' .$

HS4: If ${\rm Ext}^1_{\Lambda} \left( \Delta _{\omega} , \Delta
_{\omega '} \right) \neq 0$ then $\omega \leq \omega ' $ and it is
not true that $\omega ' \leq \omega .$

\end{defi}

\smallskip

\begin{obs} \label{alternativa} Let $\pi : \Omega \longrightarrow \Omega / \sim $ be the canonical suprajective function, $\left( \Omega / \sim , \leq \right)$ the induced partially ordered set and $\left( \Omega / \sim , \preceq \right)$ a totally ordered set where $\preceq$ extends $\leq ,$ as in Lemma \ref{growing}.

By Lemma \ref{related} and axiom HS4 if $\overline{ \omega ' } \preceq \overline{ \omega }$ then ${\rm Ext}^1_{\Lambda} \left( \Delta _{\omega} , \Delta
_{\omega '} \right) = 0,$ and by same lemma and axiom HS3 if $ \overline{ \omega ' } \prec \overline{ \omega }$ (i.e. $\overline{ \omega ' } \preceq \overline{ \omega }$ and $\overline{ \omega ' } \neq \overline{ \omega } $) then ${\rm Hom}_{\Lambda} \left( \Delta _{\omega} , \Delta
_{\omega '} \right) = 0 .$

Moreover, by Lemma \ref{related} in the definition of b-homological system we can substitute HS3 and HS4 for the following axioms:

HS3': For any total order $\preceq$ that extends $\leq ,$ if ${\rm Hom}_{\Lambda} \left( \Delta _{\omega} , \Delta
_{\omega '} \right) \neq 0$ then $\overline{\omega} \preceq \overline{\omega '} .$

HS4': For any total order $\preceq$ that extends $\leq ,$ if ${\rm Ext}^1_{\Lambda} \left( \Delta _{\omega} , \Delta
_{\omega '} \right) \neq 0$ then $\overline{\omega} \prec \overline{\omega '} .$

\end{obs}

\smallskip

\begin{defi} Given the set $\Delta = \left\{ \Delta _{\omega} \right\}_{\omega \in
\Omega},$ as in HS2, we denote by $\widetilde{\cal{F}} \left( \Delta \right)$
the full subcategory of $\Lambda -$Mod of those $M$ having a $\widetilde{\Delta}
-$filtration; we mean for this a finite sequence of submodules $$\{ 0 \} = M_0
\subsetneq M_1 \subsetneq \cdots \subsetneq M_s = M$$ such that
$M_{i} / M_{i-1} \cong \bigoplus _{\omega \in \Omega _i} \Delta _{\omega}^{\left( I_{\omega} \right)},$ where $\Omega _i$ is a finite subset of $\Omega$ and $I_{\omega}$ is an arbitray set of indexes.

The modules in $\widetilde{\cal{F}} \left( \Delta \right)$ will be called $\Delta -$filtered extended modules.

A \emph{slim} $\widetilde{\Delta}-$filtration of $M$ is a finite sequence of submodules $$\{ 0 \} = H_0
\subsetneq H_1 \subsetneq \cdots \subsetneq H_t = M$$ such that
$H_{j} / H_{j-1} \cong \Delta _{\omega _j}^{\left( I_{\omega _j} \right)},$ where $\omega _j \in \Omega$ and $I_{\omega _j}$ is an arbitrary set of indexes.

For a given $\widetilde{\Delta}
-$filtration $F$ of $M$ we denote by $\ell _{\omega} \left( F \right) = \left|  I _{\omega _{j_1} }\right| + \cdots +  \left|  I _{\omega _{j_a} }\right|,$ where $\left\{ \omega _{j_1}, \ldots , \omega _{j_a} \right\}$ are all the $\omega _j$ equal to $\omega .$ 
\end{defi}

\smallskip

\begin{obs} \label{obvio_pero_necesario}

It is immediate that $\widetilde{\cal{F}} \left( \Delta \right)$ is closed under
isomorphisms, extensions, and contains any zero object. 

Applying the Correspondence Theorem for modules it is easy to verify that for a $\widetilde{\Delta}
-$filtration $F$ of $M \in \widetilde{\cal{F}} \left( \Delta \right)$ it is possible to obtain through refinement a slim $\widetilde{\Delta}
-$filtration $\underline{F},$ and clearly $\ell _{\omega} \left( F \right)  = \ell _{\omega} \left( \underline{F} \right)$ for all $\omega \in \Omega .$

\end{obs}

\smallskip

\begin{ejem} \label{ejemplito_1} Let $\Lambda$ be an Artin $C-$algebra and $\Delta = \left\{ P_1, P_2, \ldots , P_m \right\}$ a complete list, without repetitions, of representants of isomorphism classes of indecomposable projective $\Lambda -$modules. Let $\Omega = \left\{ 1, 2, \ldots , m \right\}$ and let $i$ be related to $j$ if there is a non-zero morphism $f : P_i \longrightarrow P_j;$ the transitive closure of this relation is a pre-order $\leq$ on $\Omega .$ It is easy to verify that $\left( \Delta ; \Omega , \leq \right)$ is a b-homological system, and of course $\widetilde{\cal{F}} \left( \Delta \right)$ is the subcategory of $\Lambda -$Mod of the projectives modules.
\end{ejem}

\smallskip

\begin{ejem} \label{ejemplito_2} Let $\Lambda$ be a quasi-hereditary algebra and let $\Delta = \left\{ S_1, S_2, \ldots , S_m \right\}$ be a complete list, without repetitions, of representants of isomorphism classes of simple $\Lambda -$modules. It is known that there is a total order $\leq$ on $\Omega = \left\{ 1, 2, \ldots , m \right\}$ such that $\left( \Delta ; \Omega , \leq \right)$ is a b-homological system.
\end{ejem}

\smallskip

\begin{ejem} \label{ejemplazo} In \cite{Tr} there are examples of b-homological systems with $\Delta$ infinite: see the definition of page 16 and Theorem 1 of page 17.
\end{ejem}

\medskip

In the rest of the section we assume a given $\left( \Delta ; \Omega , \leq
\right)$ b-homological system and we will use the notation of Remark \ref{alternativa}.

\smallskip

\begin{defi} Let $\underline{F} = \left\{ H_0 \subsetneq H_1 \subsetneq \cdots
\subsetneq H_t \right\}$ be a slim $\widetilde{\Delta} -$filtration of $M \in
\widetilde{\cal{F}} \left( \Delta \right).$ The
\emph{order vector of $\underline{F}$} is defined as $h \left( \underline{F} \right) = \left( \overline{\omega}_1,\overline{\omega}_2, \ldots , \overline{\omega}_t \right),$ where $H_{j} / H_{j-1} \cong \Delta _{\omega _j}^{\left( I_{\omega _j} \right)}.$
\end{defi}

\smallskip

\begin{lema} \label{ejercicio_AH} Let $A_{i},$ with $i \in I,$ and $B_{j},$ with $j \in J,$ be $R-$modules such that $A_i$ is f.g. for any $i.$ Then for arbitrary sets of indexes $S_i$ and $T_j$ and $n \in \mathbb{N} \cup \{ 0 \}$ there is a canonical isomorphism:

$${\rm Ext}^n_{\Lambda} \left( \bigoplus _{i \in I} A_i ^{\left( S_i \right)}, \bigoplus _{j \in J} B_j^{\left( T_j \right)} \right) \cong 
\prod _{i \in I} \left[ \bigoplus _{j \in J} \left( {\rm Ext}^n_{\Lambda} \left( A_i , B_j \right) \right)^{\left( T_j \right)} \right]^{S_i} $$

In particular, if ${\rm Ext}^n_{\Lambda} \left( A_i , B_j \right) = 0$ for each $i$ and each $j$ then 
$${\rm Ext}^n_{\Lambda} \left( \bigoplus _{i \in I} A_i ^{\left( S_i \right)}, \bigoplus _{j \in J} B_j^{\left( T_j \right)} \right) = 0$$

\end{lema}

{\bf Proof:} By Theorem 7.13 of \cite{Rot} there are isomorphisms
$${\rm Ext}^n_{\Lambda} \left( \bigoplus _{i \in I} A_i ^{\left( S_i \right)}, \bigoplus _{j \in J} B_j^{\left( T_j \right)} \right) \cong 
\prod _{i \in I} {\rm Ext}^n_{\Lambda} \left( A_i ^{\left( S_i \right)}, \bigoplus _{j \in J} B_j^{\left( T_j \right)} \right) \cong 
\prod _{i \in I}\left[ \left( {\rm Ext}^n_{\Lambda} \left( A_i , \bigoplus _{j \in J} B_j^{\left( T_j \right)}  \right) \right) \right]^{S_i} $$

For $i \in I$ it is assumed that $A_i$ is f.g., an so it is known the existence of the following isomorphisms
$${\rm Ext}^n_{\Lambda} \left( A_i , \bigoplus _{j \in J} B_j^{\left( T_j \right)}  \right) \cong 
\bigoplus _{j \in J} {\rm Ext}^n_{\Lambda} \left( A_i ,  B_j^{\left( T_j \right)}  \right)  \cong 
\bigoplus _{j \in J} {\rm Ext}^n_{\Lambda} \left( A_i ,  B_j  \right)^{\left( T_j \right)}  $$

$\hfill \square$

\smallskip

\begin{prop} \label{arrangement} Let $M \in \widetilde{\cal{F}} \left( \Delta
\right) \setminus \{ 0 \}$ and $\underline{F} = \left\{ H_0 \subsetneq H_1
\subsetneq \cdots \subsetneq H_t \right\}$ be a slim $\widetilde{\Delta}
-$filtration of $M.$ 

Then there exists a slim $\widetilde{\Delta} -$filtration $\underline{F}' = \left\{
H'_0 \subsetneq H'_1 \subsetneq \cdots \subsetneq H'_t \right\}$ of
$M$ with $\ell _{\omega} \left( \underline{F} \right) = \ell _{\omega} \left( \underline{F}'
\right),$ for each $\omega \in \Omega ,$ and such that $h \left( \underline{F}'
\right) = \left( \overline{\omega '_1}, \overline{\omega '_2}, \ldots , \overline{\omega ' _t} \right)$ satisfies $\overline{\omega '_t} \preceq 
\overline{\omega '_{t-1}} \preceq \cdots \preceq \overline{\omega '_1} .$
\end{prop}

{\bf Proof:} Let us assume for some $j_0 \in \left\{ 1, 2, \ldots , t-1 \right\}$
that $\overline{\omega}_{j_0} \prec \overline{\omega}_{j_0 + 1}.$ 

Then, from HS4, Remark \ref{alternativa} and Lemma \ref{ejercicio_AH} the exact sequence
$$\xymatrix{0 \ar[r] & H_{j_0} / H_{j_0 - 1} \ar[r] & H_{j_0 + 1} / H_{j_0 - 1} \ar[r]
& H_{j_0 + 1} / H_{j_0} \ar[r] & 0}$$ splits, so there exists a submodule
$\overline{N}$ of $H_{j_0 + 1} / H_{j_0 - 1}$ such that $\overline{N} \cong
H_{j_0 + 1} / H_{j_0}$ and $\left( H_{j_0 + 1} / H_{j_0 - 1} \right) / \overline{N}
\cong H_{j_0} / H_{j_0 - 1}.$

Let $p : H_{j_0+1} \longrightarrow H_{j_0 + 1} / H_{j_0 - 1}$ be the canonical
epimorphism and $N$ the pre-image of $\overline{N}$ under $p.$ 

The following is a slim $\widetilde{\Delta}-$filtration of $M:$

$\underline{F}_1 = \left\{ H_0 \subsetneq H_1 \subsetneq
\cdots \subsetneq H_{j_0 - 1} \subsetneq N \subsetneq H_{j_0 + 1} \subsetneq
\cdots \subsetneq H_t \right\}$ 

It is easy to verify that $\ell _{\omega} \left( \underline{F}
\right) = \ell _{\omega} \left( \underline{F}_1 \right)$ for each $\omega,$ and
its order vector is equal to 
$$\left( \overline{\omega}_{1}, \ldots , \overline{\omega}_{j_0 - 1}, \overline{\omega}_{j_0 + 1},
\overline{\omega}_{j_0}, \overline{\omega}_{j_0 + 2} , \ldots , \overline{\omega}_{t} \right)$$

We can repeat this process a finite number of steps in order to reach a slim
$\widetilde{\Delta}-$filtration with a descending order vector. $\hfill \square$

\smallskip

The Proposition \ref{arrangement} provides another $\widetilde{\Delta}-$filtration; let $M \in \widetilde{\cal{F}} \left( \Delta
\right) \setminus \{ 0 \}$ and $\underline{F} = \left\{ H_0 \subsetneq H_1
\subsetneq \cdots \subsetneq H_t \right\}$ be a slim $\widetilde{\Delta}
-$filtration of $M$ such that $h \left( \underline{F}
\right) = \left( \overline{\omega _1}, \overline{\omega _2}, \ldots , \overline{\omega _t} \right)$ satisfies $\overline{\omega _t} \preceq 
\overline{\omega _{t-1}} \preceq \cdots \preceq \overline{\omega _1} .$

Assume $\overline{\omega _{j_0}} = \overline{\omega _{j_0+1}}$ then we have an exact sequence 
$$\xymatrix{0 \ar[r] & H_{j_0} / H_{j_0 - 1} \ar[r] & H_{j_0 + 1} / H_{j_0 - 1} \ar[r]
& H_{j_0 + 1} / H_{j_0} \ar[r] & 0}$$ 
and as before, by HS4, Remark \ref{alternativa} and Lemma \ref{ejercicio_AH} this exact sequence splits, so we get 
$H_{j_0 + 1} / H_{j_0 - 1} \cong \Delta _{\omega _{j_0}}^{\left( I_{\omega_{j_0}} \right)} \oplus \Delta _{\omega _{j_0 + 1}}^{\left( I_{\omega _{j_0 + 1}} \right)} .$

We can apply this idea to obtain a $\widetilde{\Delta} -$filtration
$$\{ 0 \} = W_{u_{a+1}} \subsetneq W_{u_a} \subsetneq \cdots \subsetneq W_{u_1} = M$$
where $u_{s} \in \Omega / \sim $ for $s \in \left\{ a, a-1, \ldots , 2, 1 \right\},$ $u_{1} \prec u_{2} \prec \cdots \prec u_{a-1} \prec u_{a},$  and (notice $u_{a+1}$ is just a tag) 
$W_{u_{s}} / W_{u_{s+1}} \cong \displaystyle \bigoplus _{\substack{\omega \in \Omega \\ \overline{\omega} = u_s}} \Delta _{\omega} ^{\left( I_{\omega}\right)} ,$ and in this sum for almost all index set is true $I_{\omega} = \emptyset ,$ but at least one is not empty.

\begin{defi} \label{name_only_one} A $\widetilde{\Delta} -$filtration as above will be called an \emph{ordered} $\widetilde{\Delta}-$filtration.
\end{defi}

Let $W$ be the ordered $\widetilde{\Delta}-$filtration of above; it is clear that $\ell  _{\omega} \left( W \right) = \ell _{\omega} \left( \underline{F} \right).$

\smallskip

\begin{lema} \label{restricted} Let $\{ 0 \} = W_{u_{a+1}} \subsetneq W_{u_a} \subsetneq \cdots \subsetneq W_{u_1} = M$ and  $\{ 0 \} = W_{u'_{b+1}} \subsetneq W'_{u'_b} \subsetneq \cdots \subsetneq W'_{u'_1} = N$ be ordered $\widetilde{\Delta}-$filtrations. Let $f : W_{u_a} \longrightarrow N$ be a morphism. If $u'_b \prec u_a$ then ${\rm Im} f = \{ 0 \}.$ If ${\rm Im} f \neq \{ 0 \}$ then ${\rm Im} f \subseteq W'_{u'_c},$ where $u'_c$ is the first element of $\left\{ u'_{1}, u'_{2}, \ldots , u'_{b-1}, u'_{b} \right\}$ such that $u_a \preceq u'_c .$ 
\end{lema}

{\bf Proof:} Let $\pi _1 : W'_{u'_1} \longrightarrow N / W'_{u'_2}$ be the canonical epimorphism. By HS3, Remark \ref{alternativa} and Lemma \ref{ejercicio_AH}, if $u'_1 \prec u_a$ then $\pi _1 f = 0,$ and so ${\rm Im} f \subseteq W'_{u'_2}.$ The statement can be proved repeating the argument.
$\hfill \square$

\smallskip

\begin{prop} \label{slices} Let $M \in \widetilde{\cal{F}} \left(
\Delta \right) \setminus \{ 0 \}.$ There exists an unique ordered $\widetilde{\Delta}-$filtration of $M.$ 
\end{prop}

{\bf Proof:} The existence of at least one ordered $\widetilde{\Delta}-$filtration was established in the lines before Definition \ref{name_only_one}. 

Let $\{ 0 \} = W_{u_{a+1}}  \subsetneq W_{u_a} \subsetneq \cdots \subsetneq W_{u_1} = M$ and $\{ 0 \} = W_{u_{a+1}} \subsetneq W'_{u'_{b}} \subsetneq
\cdots \subsetneq W'_{u'_1} = M$ be ordered $\widetilde{\Delta}-$filtrations.

Consider the canonical inclusion $\iota : W_{u_a} \longrightarrow M.$ By Lemma \ref{restricted} we get $u_a \preceq u'_b.$
In a similar way $u'_b \preceq u_a$ and so $u_a = u'_b.$ 

Then, applying again Lemma \ref{restricted} to both canonical inclusions, we get $W_{u_a} = W'_{u'_b}.$

We can repeat previous arguments in $M / W_{u_a}$ in order to show $u_{a-1} = u'_{b-1}$ and $W_{u_{a-1}} / W_{u_a} =
W'_{u'_{b-1}} / W'_{u'_{b}},$ so $W_{u_{a-1}} = W'_{u'_{b-1}}.$

Inductively it follows $a = b,$ $u_{j} = u'_{j}$ and $W_{u_j} = W'_{u_j},$ for $j \in \left\{ a, a-1, \ldots , 2, 1 \right\}.$ 
$\hfill \square$

\smallskip

\begin{cor} \label{well_defined} $\,$

\begin{enumerate}
\item Let $M \in \widetilde{\cal{F}} \left( \Delta \right)$ and $F$ and
$F'$ be $\widetilde{\Delta} -$filtrations of $M.$ Then $\ell _{\omega} \left( F
\right) = \ell _{\omega} \left( F' \right)$ for each $\omega \in
\Omega .$ It follows that the cardinality of factors isomorphic to $\Delta _{\omega}$ are well defined for $M,$ i.e. $\ell _{\omega} (M)$ makes sense.
\item Let $L, N \in \widetilde{\cal{F}} \left( \Delta \right)$ and $0 \rightarrow L \rightarrow M \rightarrow N
\rightarrow 0$ be an exact sequence. Then $\ell _{\omega} \left( M
\right) = \ell _{\omega} \left( L \right) + \ell _{\omega} \left( N
\right)$ for each $\omega \in \Omega .$
\end{enumerate}

\end{cor}

{\bf Proof:}

1.- As pointed out in Remark \ref{obvio_pero_necesario} $F$ induces a slim $\widetilde{\Delta}
-$filtration $\underline{F}$ and $F'$ induces a slim $\widetilde{\Delta}
-$filtration $\underline{F'} ,$ and by Proposition \ref{slices} they both induce the same ordered $\widetilde{\Delta}-$filtration $W$ of $M.$

Then, by that remark and the way to obtain $W$ we have 
$$\ell_{\omega} \left( F \right) = \ell_{\omega} \left( \underline{F} \right) = \ell_{\omega} \left( Z \right) = \ell_{\omega} \left( \underline{F'} \right) = \ell_{\omega} \left( F' \right)$$
for any $\omega.$

2.- It is a direct consequence of item (1).

$\hfill \square$

\smallskip

\begin{lema} \label{la_suma} Let $M = \bigoplus _{i \in I} G_i$ where $G_i$ is f.g and has local ring of endomorphisms for each $i.$
 Assume $M = L \oplus N.$ Then there exist decompositions in direct sums $L = \bigoplus _{j \in I'} L_j$ and $\bigoplus _{t \in I''} N_t$ such that:

\begin{enumerate} 

\item For any $j$ there exists an $i_j \in I$ such that $L_j \cong G_{i_j}.$
 
\item For any $t$ there exists an $i_t \in I$ such that $N_t \cong G_{i_t}.$

\item Let $G$ be an $R-$module with local ring of endomorphisms and $\ell _G (M)$ the cardinal of the subset $I_G$ of $I,$ where $i \in I_G$ iff $G_i \cong G.$ Then $\ell _G (M) = \ell _G (L) + \ell _G (N).$
 
\end{enumerate}

\end{lema}

{\bf Proof:} By the Crawley-J$\o$nsson-Warfield Theorem (26.5 of \cite{AF}) and a review of its proof, we get items (1) and (2).

The third item follows by Azumaya's Decomposition Theorem (12.6 of \cite{AF}).
$\hfill \square$

\smallskip

For the next proof recall that given
$L, N \in R -$Mod, the \emph{trace} ${\rm tr}_{L} \left( N
\right)$ of $L$ on $N$ is the sum of all the images of the
homomorphisms from $L$ to $N.$

\smallskip

\begin{cor} \label{closed}
$\widetilde{\cal{F}} \left( \Delta \right)$ is closed under direct summands.
\end{cor}

{\bf Proof:}  Let $M \in \widetilde{\cal{F}} \left( \Delta \right)$ and $M = L \oplus N.$

Consider the ordered $\widetilde{\Delta}-$filtration of $M$ with the notation of Proposition \ref{slices}.

By Lemma \ref{restricted} we have the identity $W_a = {\rm tr}_{W_a}
\left( M \right),$ and by additivity of the trace, ${\rm tr}_{W_a}
\left( M \right) = {\rm tr}_{W_a} \left( L \right) \oplus {\rm
tr}_{W_a} \left( N \right).$ 

By Lemma \ref{la_suma} we get that ${\rm tr}_{W_a} \left( L
\right) = L_a \subseteq L$ and ${\rm tr}_{W_a} \left( N
\right) = N_a \subseteq N$ are in $\widetilde{\cal{F}} \left(
\Delta \right) .$

By same Lemma $\ell _{\omega} \left( W_a \right) = \ell _{\omega} \left( L_a \right) + \ell _{\omega} \left( N_a \right)$ for each $\omega .$

We can repeat this argument for the quotients $M /W_a \cong \left(
L / L_a \right) \oplus \left( N / N_a \right)$ in order to get larger
submodules of $L$ and $N$ that are in $\widehat{\cal{F}} \left(
\Delta \right).$

Repeating the procedure we obtain ordered $\widetilde{\Delta}-$filtrations of $L$ and
$N,$ so they both belong to $\widetilde{\cal{F}} \left( \Delta \right).$
$\hfill \square$

\end{document}